\documentclass{article}
\usepackage{amssymb,amsthm}
\usepackage{stmaryrd}
\usepackage{mathtools}
\usepackage{amsfonts}
\usepackage{amsmath}
\usepackage{txfonts,dsfont}
\DeclareMathAlphabet{\mathcal}{OMS}{cmsy}{m}{n} 
\usepackage{enumerate}
\usepackage{comment}
\usepackage[colorlinks=true, allcolors=blue]{hyperref}
\allowdisplaybreaks[4]
\usepackage{color}
\usepackage{pdfsync}
\usepackage[numbers,sort&compress]{natbib}
\usepackage[letterpaper,top=2cm,bottom=2cm,left=3cm,right=3cm,marginparwidth=1.75cm]{geometry}
\usepackage{graphicx}
\usepackage{titlesec}
\titleformat{\section}{\normalfont\Large\bfseries\raggedright}{\thesection.}{0.5em}{}
\titleformat{\subsection}{\normalfont\large\bfseries\raggedright}{\thesubsection.}{0.5em}{}

\newtheorem{theorem}{Theorem}[section]
\newtheorem{lemma}[theorem]{Lemma}
\newtheorem{proposition}[theorem]{Proposition}
\newtheorem{corollary}[theorem]{Corollary}
\newtheorem{conjecture}[theorem]{Conjecture}
\theoremstyle{definition}

\theoremstyle{remark}
\newtheorem{remark}[theorem]{Remark}
\numberwithin{equation}{section}
\newcommand{\R}{\mathbb R}

\newcommand{\W}{\mathcal W}
\newcommand{\dd}{\,\mathrm d}
\newcommand{\Id}{\mathrm{Id}}
\newcommand*{\Rmn}[1]{\uppercase\expandafter{\romannumeral#1}}
\DeclareMathOperator{\diver}{div}

\DeclareMathOperator{\tr}{tr}
\DeclareMathOperator{\Ric}{Ric}

\makeatletter
\@namedef{subjclassname@2020}{\textup{2020} Mathematics Subject Classification}
\makeatother

\begin{document}

\title{Counterexamples to Ecker's Conjecture}

\author{
Jiarui Chen \footnote{1. College of Mathematics and Statistics, Northwest Normal University, Lanzhou, 730070, Gansu, China; 2. Gansu Provincial Research Center for Basic Disciplines of Mathematics and Statistics, Lanzhou, 730070, Gansu, China. This work is supported by the Fund for Innovative Fundamental Research Group Project of Gansu Province (No. 24JRRA778). \texttt{jiaruichen@nwnu.edu.cn}}  \quad and \quad Jie Wang \footnote{School of Mathematics and Statistics, Jiangsu Normal University, Xuzhou, 221116, Jiangsu, China. (Corresponding author). \texttt{jiewang@jsnu.edu.cn}}
}

\date{}
	
\maketitle

\vskip2pt\vskip2pt

\noindent
{\bf Abstract.}
Analogous to Perelman's $\mathcal{W}$-entropy defined for Ricci flow, K. Ecker in \cite{Ecker07} defined a functional $\mathcal{W}$ for the mean curvature flow in $\mathbb{R}^{n+1}$ and the region it encloses, and made the conjecture that this functional is monotonically increasing in time. In this paper, we give a negative answer to Ecker’s conjecture.
\vskip2pt\vskip2pt



\vskip2pt\vskip2pt

\section{Introduction}

Let $(M,g)$ be a closed $n$-dimensional Riemannian manifold. Perelman in his 2002 preprint~\cite{Perelman} introduced the well-known $\mathcal{W}$-entropy functional
\[
 \W(g,f,\tau):=\int_M\bigl[\tau(|\nabla f|^2+R_g)+f-n\bigr]u\dd V_g,
 \qquad u=(4\pi\tau)^{-n/2}e^{-f},
\]
where $f\in C^{\infty}(M)$, $\tau>0$ is a parameter, $R_g$ denotes scalar curvature and $\int_M u\dd V_g=1$. 

For the Ricci flow $\partial_tg=-2\Ric_g$ on $M$ with 
\[
 \partial_t f=-\Delta_gf+|\nabla f|^2-R_g+\frac {n}{2\tau},\qquad \frac{d}{dt}\tau=-1,
\]
Perelman proved the following monotonicity formula holds
\[
 \frac{\mathrm d}{\mathrm dt}\W
 =2\tau\int_M\left|\Ric_g+\nabla^2f-\frac g{2\tau}\right|^2u\dd V_g
 \geq0.
\]
This formula turns out being of fundamental importance in the study of Ricci flow.

Since then, entropy formulas have played an important role in the study of various geometric evolution equations. For the linear heat equation $\partial_t u=\Delta_gu$ on a fixed manifold $(M^n,g)$ with $u=(4\pi \tau)^{-n/2}e^{-f}$, $\int_M u(\cdot,0)\dd V_g=1$ and $d\tau/dt=1$, Ni~\cite{Ni04} obtained a corresponding formula 
\[
 \frac{\mathrm d}{\mathrm dt}\int_M(\tau|\nabla f|^2+f-n)u\dd V_g
 =-2\int_M \tau\left(\left|\nabla^2f-\frac g{2\tau}\right|^2
                   +\Ric_g(\nabla f,\nabla f)\right)u\dd V_g.
\]
In particular, nonnegative Ricci curvature gives a nonincreasing entropy. Li~\cite{Li12} extended this result to the Witten Laplacian under Bakry--\'Emery Ricci curvature assumptions. These results suggest looking for similar formulas for other geometric flows.

In \cite{Ecker07}, Ecker found a surprising link between Perelman's ${\mathcal W}$-functional for Ricci flow and Hamilton's  \cite{Hamilton95} differential Harnack expression for mean curvature flow in $\R^m$; he considered instead an entropy of the enclosed region, with a boundary term involving its normal speed. In the following, we recall Ecker's definition in the case of mean curvature flow more precisely.

Throughout the paper, $n\geq1$ is the dimension of a hypersurface  and $m=n+1$ is the dimension of ambient spaces.
Let $\left(\Omega_t\right)_{t\in [0,T)}\subset\R^m$ be a family of bounded open set with smooth compact
boundary $M_t=\partial \Omega_t$. Let $\nu$ be the outward unit normal, the second fundamental form $A$ and mean curvature $H$ are defined by
\[
 A(V,W):=\langle D_V\nu,W\rangle,\qquad H:=\tr A.
\]
By this convetion, $H=n/R$ on a round sphere of radius $R$, and the normal mean curvature
flow is $\partial_tX=-H\nu$. 

For $\tau=a-t>0$, let
\begin{equation}\label{eq:normalization}
 u=(4\pi\tau)^{-m/2}e^{-f},\qquad \int_{\Omega_t}u\dd x=1.
\end{equation}
Assume that $f$ satisfies the following evolution equation in $\Omega_t$  with Neumann boundary condition
\begin{equation}\label{eq:system}
 \begin{cases}
 \partial_t f+\Delta f=|\nabla f|^2+\dfrac m{2\tau},&\text{in }\Omega_t,\\[2mm]
 \partial_\nu f=H,&\text{on }M_t.
 \end{cases}
\end{equation}
Let $\Omega_t$ evolve by diffeomorphisms generated by $-\nabla f$. Its normal component on $M_t$ is $-H\nu$, so the boundary hypersurfaces $M_t$ are just evolving along the mean curvature flow. Ecker \cite{Ecker07} defined the entropy as
\begin{equation}\label{eq:W}
 \W_H(\Omega,f,\tau)
 :=\int_\Omega(\tau|\nabla f|^2+f-m)u\dd x
       +2\tau\int_{\partial\Omega}Hu\dd S,
\end{equation}
and proved that 
\begin{equation}\label{eck}
\begin{aligned}
\frac{\mathrm d}{\mathrm dt}\W_H(\Omega_t,f(t),\tau(t))&=2\tau\int_{\Omega_t}\left|\nabla^2f-\frac{\Id}{2\tau}\right|^2u\dd x\\
       &\quad\quad+2\tau\int_{M_t}\left(\frac{\partial H}{\partial t}-2\langle\nabla^MH,\nabla^Mf\rangle+A(\nabla^Mf,\nabla^Mf)-\frac H{2\tau}\right)u\dd S,
\end{aligned}
\end{equation}
where $\nabla^M$ denotes the covariant derivative on the hypersurface $M_t$. 

Notice that the integrand of the boundary term in \eqref{eck} is closely related to Hamilton's differerential Harnack expression for mean curvature flow \cite{Hamilton95}. Indeed, Hamilton's Harnack inequality states that any weakly convex solution $X_t:M\rightarrow \mathbb{R}^{m}$ of the mean curvature flow 
\[
\partial_t X=-H\nu,
\] 
satisfies that for all $t\in [0,T)$ and all tangent vectors $V$ on $M_t=X_t(M)$, the following inequality holds:
\[
    \frac{\partial H}{\partial t}+2\langle\nabla^MH,V\rangle+A(V,V)+\frac{H}{2t}\geq 0,
\]
where $A$ is the second fundamental form of $M_t$. If we denote by
\[
Z(V):=\frac{\partial H}{\partial t}+2\langle\nabla^MH,V\rangle+A(V,V),
\]
then Ecker's evolution formula \eqref{eck} can be written as
\[
\begin{aligned}
\frac{\mathrm d}{\mathrm dt}\W_H(\Omega_t,f(t),\tau(t))&=2\tau\int_{\Omega_t}\left|\nabla^2f-\frac{\Id}{2\tau}\right|^2u\dd x+2\tau\int_{M_t}\left(Z(-\nabla^Mf)-\frac H{2\tau}\right)u\dd S.
\end{aligned}
\]

One can observe that the first term is the integral over $\Omega_t$ of a nonnegative quantity, as in Perelman's work. The second term is an integral related to Hamilton's differerential Harnack expression over $\partial \Omega_t$. Motivated by this relation, Ecker proposed the following conjecture.

\begin{conjecture}[Ecker, \cite{Ecker07}]\label{conj:ecker}
For a mean curvature flow of compact embedded hypersurfaces $M_t$ in $\R^{n+1}$ and the above $f(t),\tau(t)$, the inequality
\begin{equation}\label{eck2}
\begin{aligned}
\frac{\mathrm d}{\mathrm dt}\W_H(\Omega_t,f(t),\tau(t))&=2\tau\int_{\Omega_t}\left|\nabla^2f-\frac{\Id}{2\tau}\right|^2u\dd x+2\tau\int_{M_t}\left(Z(-\nabla^Mf)-\frac H{2\tau}\right)u\dd S\geq 0
\end{aligned}
\end{equation}
holds for all $t\in[0,T)$.
\end{conjecture}


Ecker's later work~\cite{Ecker14} further develops entropy and differential Harnack formulas for evolving domains, including estimates that allow controlled negative entropy derivatives. Lott \cite{Lott12} established analogous relations for mean curvature flow in an arbitrary Ricci flow background. In flat case, Lott \cite{Lott12} obtained the ${\mathcal F}$-entropy functional counterpart of Ecker's result, which can date back to Perelman's work on ${\mathcal F}$-functional for Ricci flow \cite{Perelman}. Recently, by modifying Ecker's definition and using Hamilton's Harnack inequality for mean curvature flow, Li and Yan~\cite{LiYan25} proved that their redefined $\mathcal{W}$-entropy is monotonically decreasing in time. Additionally, they provide a rigidity theorem for their $\mathcal{W}$-entropy.

In this paper, we give a negative answer to Ecker’s conjecture by an explicit calculation on a family of shrinking round balls. 

\begin{theorem}\label{thm:main}
Assume $S>0$ and $a>0$. Let $s=S-t$ and $\tau=a-t$. On the balls
$\Omega_t=B_{\sqrt{2ns}}(0)\subset\R^{n+1}$, set
\begin{equation}\label{eq:main-data}
 f(x,t)=\frac{|x|^2}{4s}+\frac {n+1}{2}\log\frac s\tau+\log P_n,\qquad u(x,t)=\frac{e^{-|x|^2/(4s)}}{P_n(4\pi s)^{\frac{n+1}{2}}},
\end{equation}
where 
\[
P_n:=\int_{B_{\sqrt{2n}}(0)}(4\pi)^{-\frac{n+1}{2}}e^{-|y|^2/4}\dd y.
\]
Then for $t\in [0,\min\{S,a\})$, $\Omega_t, f(x,t), u(x,t),\tau$ satisfy all the assumptions of Conjecture \ref{conj:ecker}, and
\[
\frac{\mathrm d}{\mathrm dt}\W_H(\Omega_t,f(t),\tau(t))=2\tau\left((n+1)\left(1-\frac{\tau}{s}\right)-2\kappa_n \frac{\tau}{s}\right)\cdot\frac{1-\frac{\tau}{s}}{4s^2\left(\frac{\tau}{s}\right)^2},
\]
where
\[
\kappa_n:=\frac{|\mathbb S^n|\,n(2n)^{(n-1)/2}e^{-n/2}}
                    {(4\pi)^{\frac{n+1}{2}}P_n}.
\]
Consequently, 
\[
\frac{\mathrm d}{\mathrm dt}\W_H(\Omega_t,f(t),\tau(t))<0\quad\quad \text{if and only if}\quad\quad \frac {n+1}{n+1+2\kappa_n}<\frac{\tau}{s}<1.
\]
Hence, on every interval on which this condition holds,
$\W_H\left(\Omega_t,f(t),\tau(t)\right)$ is strictly decreasing.
\end{theorem}

\begin{corollary}
For every $n\geq1$, choose $S=1$ and $a=\frac{n+1}{n+2}$. Then for the examples in Theorem \ref{thm:main}, we have
\[
 \left.\frac{\mathrm d}{\mathrm dt}\right|_{t=0}\W_H(\Omega_t,f(t),\tau(t))=2\tau\cdot\frac{1-2\kappa_n}{4(n+1)}<0,\qquad
 \frac{\mathrm d}{\mathrm dt}\W_H(\Omega_t,f(t),\tau(t))<0\quad\text{for }0\leq t<t^*,
\]
where
\[
 t^* =\frac{(n+1)(2\kappa_n-1)}{2\kappa_n(n+2)},\qquad 0<t^*<a.
\]
\end{corollary}

We give an intuitive explanation for the construction of this counterexample. The boundary condition fixes the parameter $s$, but it leaves the entropy parameter $\tau$ free. When $\tau$ is slightly smaller than $s$, the entropy parameter vanishes before the sphere becomes singular and the boundary contribution to~\eqref{eck2} is negative to first order in $s-\tau$, whereas the interior term is quadratic. This is the source of the counterexample.

\begin{theorem}\label{thm:first-singular}
For every $n\geq1$, there is a two-component mean curvature flow on its maximal smooth interval $[0,T)$, with $\tau=T-t$ and $\Omega_t, f(x,t), u(x,t),\tau$ satisfy all the assumptions of Conjecture \ref{conj:ecker}, such that
\[
\left.\frac{\mathrm d}{\mathrm dt}\right|_{t=0}\W_H(\Omega_t,f(t),\tau(t))<0,
\]
and 
\[
\frac{\mathrm d}{\mathrm dt}\W_H(\Omega_t,f(t),\tau(t))<0,\quad\text{for \;\;} 0\leq t<t^*,
\]
where $0<t^*<T$ is a constant.
\end{theorem}

The paper is organized as follows. In Section~\ref{se2}, we
give the proof of Theorem \ref{thm:main}. In Section~\ref{se3}, we prove Theorem~\ref{thm:first-singular}. In Section~\ref{se4}, we show that shrinking embedded circles exactly provide concrete counterexamples to Ecker's conjecture in $\mathbb{R}^2$.

\section{Proof of Theorem \ref{thm:main}}\label{se2}

We first reformulate~\eqref{eck2} to the following form.

\begin{proposition}\label{prop:boundary-defect}
Set $B:=\nabla^2f-\frac{\Id}{2\tau}$. Then
\begin{equation}\label{eq:boundary-defect}
 \int_{\Omega_t}\left|\nabla^2f-\frac{\Id}{2\tau}\right|^2u\dd x+\int_{M_t}\left(Z(-\nabla^Mf)-\frac H{2\tau}\right)u\dd S=\int_{\Omega_t}|B|^2u\dd x
        +\int_{M_t}\langle A,B|_{TM_t}\rangle u\dd S,
\end{equation}
where $B|_{TM_t}$ is the restriction of $B$ to tangent vectors on $M_t$.
\end{proposition}

\begin{proof}
The mean curvature evolution equation and the contracted Codazzi
identity give
\[
 H_t=\Delta_MH+|A|^2H,\qquad \diver_M A=\nabla^MH,
\]
see~\cite[Appendices~A and~B]{EckerBook}. Put $X=\nabla^Mf$.
Since $\nabla^Mu=-uX$, integration by parts on the closed hypersurface
gives
\[
 \int_{M_t}(\Delta_MH-2\langle\nabla^MH,X\rangle)u\dd S
 =-\int_{M_t}\langle\nabla^MH,X\rangle u\dd S.
\]
Similarly, integrating $\diver_M(uA(X,\cdot))$ yields
\[
 \int_{M_t}\bigl(-\langle\nabla^MH,X\rangle+A(X,X)\bigr)u\dd S
 =\int_{M_t}\langle A,\nabla_M^2f\rangle u\dd S.
\]
Since $\partial_\nu f=H$, then we have
\[
 \nabla_M^2f=(\nabla^2f)|_{TM_t}-HA.
\]
The terms involving $H|A|^2$ therefore cancel, and we obtain
\[
 \int_{M_t}Z(\nabla^Mf)u\dd S
 =\int_{M_t}\langle A,(\nabla^2f)|_{TM_t}\rangle u\dd S.
\]
Subtracting $\int_{M_t}\frac{Hu}{2\tau}\dd S$ proves the proposition.
\end{proof}

From Proposition \ref{prop:boundary-defect}, one can see that the first term in the right hand side of~\eqref{eq:boundary-defect} is quadratic in $B$, while the second is linear. In particular, positive definiteness of $A$ does not determine the sign of the latter term. 

Next, we construct the counterexamples. Let $S>0$ and $a>0$. Let $s=S-t$ and $\tau=a-t$. We consider
\begin{equation}\label{eq:gaussian}
 \begin{aligned}
 \Omega_t&=B_{\sqrt{2ns}}(0),\qquad M_t=\partial\Omega_t,\\
 u(x,t)&=\frac{e^{-|x|^2/(4s)}}{P_n(4\pi s)^{\frac{n+1}{2}}},\\
 f(x,t)&=\frac{|x|^2}{4s}+\frac{n+1}{2}\log\frac s\tau+\log P_n,
 \end{aligned}
\end{equation}
for $0\leq t<\min\{S,a\}$, where 
\[
P_n=\int_{B_{\sqrt{2n}}(0)}(4\pi)^{-\frac{n+1}{2}}e^{-|y|^2/4}\dd y.
\] 
We then check all the conditions in the conjecture and compute the expressions needed below.

\begin{lemma}\label{lem:admissible}
The hypersurfaces $M_t$ in~\eqref{eq:gaussian} evolve by mean curvature flow, and $\Omega_t, f(x,t), u(x,t),\tau$ satisfy all the assumptions of Conjecture \ref{conj:ecker}. Moreover,
\begin{equation}\label{eq:moments}
 \int_{M_t}Hu\dd S=\frac{\kappa_n}s,\qquad
 \int_{\Omega_t}\frac{|x|^2}{4s}u\dd x=\frac{n+1}{2}-\kappa_n,
\end{equation}
where
\begin{equation}\label{con1}
\kappa_n:=\frac{|\mathbb S^n|\,n(2n)^{(n-1)/2}e^{-n/2}}
                    {(4\pi)^{\frac{n+1}{2}}P_n},
\end{equation}
and $\frac1 2 <\kappa_n<\frac{n+1}{2}$.
\end{lemma}

\begin{proof}
For $R=\sqrt{2ns}$, one has $dR/dt=-n/R=-H$. The change of variables $x=\sqrt{s}\,y$ proves \eqref{eq:normalization}. Direct calculations give
\begin{equation}\label{eq:derivatives}
 \begin{gathered}
 \nabla f=\frac x{2s},\qquad \nabla^2f=\frac{\Id}{2s},\qquad
 \Delta f=\frac {n+1}{2s},\qquad |\nabla f|^2=\frac{|x|^2}{4s^2},\\
 \partial_t f=\frac{|x|^2}{4s^2}-\frac {n+1}{2s}+\frac {n+1}{2\tau}.
 \end{gathered}
\end{equation}
On the boundary, we have
\[
 \partial_\nu f=\frac R{2s}=\frac nR=H.
\]
Hence the equation in~\eqref{eq:system} holds.
The maps
\begin{equation}\label{eq:transport}
 \Phi_{t,t_0}(x)=\sqrt{\frac{s(t)}{s(t_0)}}\,x
\end{equation}
are diffeomorphisms between the corresponding closed balls, with
\[
 \partial_t\Phi_{t,t_0}=-\nabla f(\Phi_{t,t_0},t).
\]
All functions are smooth up to the spatial boundary on compact time intervals before $\min\{S,a\}$, and $u$ is positive there.

Evaluating the boundary integral explicitly gives the first identity in~\eqref{eq:moments}, with $\kappa_n$ as in~\eqref{con1}. For the second identity, we observe that
\[
 \diver(xu)=\left(n+1-\frac{|x|^2}{2s}\right)u,
 \qquad \langle x,\nu\rangle=2sH\quad\text{on }M_t.
\]
The divergence theorem therefore gives
\[
 n+1-2\int_{\Omega_t}\frac{|x|^2}{4s}u\dd x=2\kappa_n.
\]
The integrand $|x|^2/(4s)$ lies strictly between $0$ and $n/2$ almost
everywhere. Its average against the positive probability density $u$
lies between the same endpoints, which proves the stated bounds on
$\kappa_n$.
\end{proof}

\begin{proposition}\label{prop:sign}
For the construction in~\eqref{eq:gaussian}, we have
\begin{equation}\label{eq:exactQ}
 \begin{aligned}
 \frac{\mathrm d}{\mathrm dt}\W_H(\Omega_t,f(t),\tau(t))=2\tau\left((n+1)\left(1-\frac{\tau}{s}\right)-2\kappa_n \frac{\tau}{s}\right)\frac{1-\frac{\tau}{s}}{4s^2\left(\frac{\tau}{s}\right)^2}.
 \end{aligned}
\end{equation}
It follows that $\frac{\mathrm{d} \W_H}{\mathrm dt}<0$ if and only if $\frac {n+1}{n+1+2\kappa_n}<\frac{\tau}{s}<1$. 
\end{proposition}

\begin{proof}
Note that both $H$ and $f$ are constant on $M_t$, and
$\frac{\partial H}{\partial t}=n^2/R^3=H/(2s)$. The boundary term in~\eqref{eck2} is therefore
\[
 2\tau\frac12\left(\frac1s-\frac1\tau\right)\int_{M_t}Hu\dd S.
\]
Using~\eqref{eq:derivatives}, $|\Id|^2=n+1$, and the normalization \eqref{eq:normalization} gives 
\[
\begin{aligned}
 \frac{\mathrm d}{\mathrm dt}\W_H(\Omega_t,f(t),\tau(t))=2\tau \frac {n+1}{4}\left(\frac1s-\frac1\tau\right)^2
     +2\tau\frac{\kappa_n}{2s}\left(\frac1s-\frac1\tau\right),
 \end{aligned}
\]
which proves~\eqref{eq:exactQ}. The remaining part follows by factoring and
determines its sign.
\end{proof}

\begin{corollary}\label{cor:dimension}
For every $n\geq1$, choose $S=1$ and $a=\frac{n+1}{n+2}$. Then
\begin{equation}\label{eq:explicit-time}
 \left.\frac{\mathrm d}{\mathrm dt}\right|_{t=0}\W_H(\Omega_t,f(t),\tau(t))=2\tau\frac{1-2\kappa_n}{4(n+1)}<0,
 \end{equation}
 and
 \begin{equation}\label{eq:explicit-time2}
 \frac{\mathrm d}{\mathrm dt}\W_H(\Omega_t,f(t),\tau(t))<0\quad\text{for \;\;}0\leq t<t^*,
\end{equation}
where
\begin{equation}\label{eq:tstar}
 t^* =\frac{(n+1)(2\kappa_n-1)}{2\kappa_n(n+2)},\qquad 0<t^*<a.
\end{equation}
\end{corollary}

\begin{proof}
At $t=0$, one has 
\[
\frac1 s-\frac1 \tau=-\frac1{n+1}.
\]
 Along the flow,
 \[
 \frac{\tau}{s}=\frac{a-t}{1-t}
 \]
 decreases from $a$ to zero. Solving $\frac{\tau}{s}=\frac {n+1}{n+1+2\kappa_n}$ gives~\eqref{eq:tstar}, and Proposition~\ref{prop:sign} applies.
\end{proof}

\begin{remark}
For any $0<T_0<t^*$, these data are smooth on $[0,T_0]$, with
$a>T_0$ and $\frac{\mathrm{d} \W_H}{\mathrm dt}<0$ throughout. Thus they also satisfy the formulation
on a closed time interval with $\tau=a-t$ and $a\geq T_0$ used in
\cite[Proposition~3.4]{Ecker07}.
\end{remark}

\section{Proof of Theorem \ref{thm:first-singular}}\label{se3}

 We consider a disjoint union of shrinking balls. In this case, the first singular time $T$ is determined by the smallest ball and the parameter $s$ on a larger component can therefore have a geometric scale greater than $\tau=T-t$. 
 
 \begin{theorem}
For every $n\geq1$, there is a two-component mean curvature flow on its maximal smooth interval $[0,T)$, with $\tau=T-t$ and $\Omega_t, f(x,t), u(x,t),\tau$ satisfy all the assumptions of Conjecture \ref{conj:ecker}, such that
\[
\left.\frac{\mathrm d}{\mathrm dt}\right|_{t=0}\W_H(\Omega_t,f(t),\tau(t))<0,
\]
and 
\[
\frac{\mathrm d}{\mathrm dt}\W_H(\Omega_t,f(t),\tau(t))<0,\quad\text{for \;\;} 0\leq t<t^*,
\]
where $0<t^*<T$ is a constant.
\end{theorem}

\begin{proof}
Take two spheres $\partial B_{\sqrt{\frac{2n(n+1)}{n+2}}}(p_1)$ and $\partial B_{\sqrt{2n}}(p_2)$ with $p_1,p_2\in\R^{n+1}$ so that the two balls they contained have disjoint closures. Consider the disjoint union of $\partial B_{\sqrt{\frac{2n(n+1)}{n+2}}}(p_1)$ and $\partial B_{\sqrt{2n}}(p_2)$ and let the union evolve by mean curvature flow. Then the first sphere collapses at time $T=\frac{n+1}{n+2}$ and its curvature becomes unbounded, while the second still has positive radius. Hence $T=\frac{n+1}{n+2}$ is the maximal smooth time of the original two-component mean curvature flow.

Define $s_1=\tau=T-t,s_2=1-t$,
\begin{equation}
 \begin{aligned}
 u_1(x,t)&=\frac1 2 \frac{ e^{-|x-p_1|^2/(4s_1)}}{P_n(4\pi s_1)^{\frac{n+1}{2}}},\\
 f_1(x,t)&=\frac{|x-p_1|^2}{4s_1}
       +\frac{n+1}{2}\log\frac{s_1}{\tau}+\log P_n+\log 2,
 \end{aligned}
\end{equation}
and
\begin{equation}
 \begin{aligned}
 u_2(x,t)&=\frac1 2 \frac{ e^{-|x-p_2|^2/(4s_2)}}{P_n(4\pi s_2)^{\frac{n+1}{2}}},\\
 f_2(x,t)&=\frac{|x-p_2|^2}{4s_2}
       +\frac{n+1}{2}\log\frac{s_2}{\tau}+\log P_n+\log 2.
 \end{aligned}
\end{equation}

On the first component, $s_1=\tau$, so by Proposition \ref{prop:sign} we have
\begin{equation}\label{w1}
\frac{\mathrm d}{\mathrm dt}\W_1 =0 \quad \text{and}\quad \W_1=\log P_n.
\end{equation}
The second component is the example in Corollary \ref{cor:dimension}. 

Lemma~\ref{lem:admissible} applies to each translated ball, with
transport maps
\[
 x\longmapsto p_i+\sqrt{\frac{s_i(t)}{s_i(t_0)}}(x-p_i).
\]
The constant $\log 2$ changes none of the derivatives of $f_i$.
Integrating component by component gives
\begin{equation}\label{eq:mixture}
 \W_H=\frac1 2 \W_1+\frac1 2 \W_2+\log 2.
\end{equation}
Therefore, \eqref{eq:mixture} together with \eqref{w1} and Corollary \ref{cor:dimension} give the desired result.
\end{proof}

\section{A concrete example in the plane}\label{se4}

For $n=1$, all constants can be written in elementary form. Let
\begin{equation}\label{eq:circle}
 \begin{gathered}
 s=1-t,\qquad \tau=\frac34-t,\qquad
 \Omega_t=B_{\sqrt{2s}}(0)\subset\R^2,
 \qquad 0\leq t\leq\frac14.
  \end{gathered}
\end{equation}

 Direct calculations gives $P_1=1-e^{-1/2}$ and $\kappa_1=1/[2(e^{1/2}-1)]$. Let
 \begin{equation}\label{eq:circle2}
 \begin{gathered}
 u(x,t)=\frac{e^{-|x|^2/(4s)}}{4\pi s P_1 },\qquad
 f(x,t)=\frac{|x|^2}{4s}+\log\frac s\tau+\log P_1 .
 \end{gathered}
\end{equation}
 In particular,
\[
\left.\frac{\mathrm d}{\mathrm dt}\right|_{t=0}\W_H(\Omega_t,f(t),\tau(t))=\frac{3}{2}\left(\frac1{18}-\frac1{12(e^{1/2}-1)}\right)<0.
\]
On the whole interval $t\in [0,\frac{1}{4}]$, we have $\frac{2}{3}\leq \frac{\tau}{s}\leq \frac{3}{4}$ and
\[
 \frac{\mathrm d}{\mathrm dt}\W_H(\Omega_t,f(t),\tau(t))=2\tau\cdot\frac{1-\frac{\tau}{s}}{2s^2\left(\frac{\tau}{s}\right)^2}\left(1-\frac{\tau}{s}-\kappa_1\left(\frac{\tau}{s}\right)\right)<0,
\]
since $(1-\frac{\tau}{s})/\frac{\tau}{s}\leq1/2<\kappa_1$. Thus shrinking embedded circles provide counterexamples to Conjecture \ref{conj:ecker} in $\mathbb{R}^2$.

\vskip12pt
\vskip12pt

\noindent{\bf
Data availability.}   All relevant data are within the paper.

\vskip12pt\vskip12pt

\noindent{\bf
Conflict of interest.}   There are no conflicts of interest to this work.

\end{document}